# Double orbits of weakly almost periodic functions

Ching Chou

**Abstract**

For a locally compact group $G$, let $AP(G)$ and $WAP(G)$ be respectively the $C^*$-algebras of almost periodic and weakly almost periodic functions on $G$. For a bounded continuous function $f$ on G, $f$ is said to be strictly w.a.p. if its double orbit $O(f)$ is relatively weakly compact and $f$ said to be strictly uniformly continuous if its double orbit is equicontinuous on $G$. The $C^*$-algebras of such functions are denoted, respectively, by $WS(G)$ and $UCS(G)$. Then $WS(G) \subset UCS(G)$ and $AP(G) \subset WS(G) \subset WAP(G)$. $G$ is called a $WS$-group if $WS(G) = WAP(G)$. We will show that if a discrete $FC$-group $G$ is a $WS$-group, then its center is of finite index in $G$. A noncompact locally compact group $G$ is minimally w.a.p., if $WAP(G) = AP(G) \oplus C_0(G)$. If $G$ is minimally w.a.p., then $WS(G) = AP(G)$, i.e., if the double orbit of a bounded continuous function $f$ is relatively weakly compact then it is relatively norm compact. It is known that for $n \geq 2$, the motion group $M(n)$, and the special linear group $SL(n, \mathbb{R})$ are minimally w.a.p. On the other hand, there exist locally compact groups $G$ such that $WS(G) = AP(G)$ but $G$ is not minimally w.a.p. We will show that if $G$ is an $IN$-group and $K = K_G$ is the intersection of all closed invariant neighborhoods of the identity of $G$, then $UCS(G) = UCS(G/K)$ and $WS(G) = WS(G/K)$. We will identify the strictly w.a.p. functions on the $ax + b$ group. We will also show that $UCS(SL(2, \mathbb{R}))$ only contains the constant functions.

## 1. Introduction

Let $G$ be a locally compact group, $C(G)$ the $C^*$-algebra of bounded complex-valued continuous functions on $G$ with the sup norm and $C_0(G)$ the C*-subalgebra of $C(G)$ consisting of functions vanishing at infinity. For $f \in C(G)$ and $x, y \in G$, the left translation of $f$ by $x$, the right translation of $f$ by $y$, and the two-sided translation of $f$ by $x$ and $y$, are respectively defined by $_xf(z) = f(xz)$, $f_y(z) = f(zy)$ and $_xf_y(z) = f(xzy)$, $z \in G$. Let $O_L(f) = \{_xf : x \in G\}$, $O_R(f) = \{f_y : y \in G\}$, and $O(f) = \{_xf_y : x, y \in G\}$ be, respectively, the left orbit, the right orbit and the double orbit of $f \in C(G)$.

For $f \in C(G)$, it is well-known and is very easy to prove that the following three conditions are equivalent: (1) $O_L(f)$ is relatively compact in $C(G)$; (2) $O_R(f)$ is relatively compact in $C(G)$; (3) $O(f)$ is relatively compact in $C(G)$. If $f \in C(G)$ satisfies one of these three equivalent conditions, then $f$ is said to be almost periodic and the set of all such functions on $G$ is denoted by $AP(G)$. Then $AP(G)$ is a $C^*$-subalgebra of $C(G)$. It is a well-known result of von Neumann that the linear span of the coefficient functions of finite dimensional continuous irreducible unitary representations of $G$ is uniformly dense in AP(G); see von Neumann [28].

By the Grothendieck weak compactness criterion [16], for $f \in C(G)$, the following two conditions are equivalent: (1)' $O_L(f)$ is relatively weakly compact in $C(G)$; (2)' $O_R(f)$ is relatively weakly compact in $C(G)$. If $f \in C(G)$, then $f$ is said to be weakly almost periodic (w.a.p.), if it satisfies (1)', or, equivalently, (2)'. The space of all continuous w.a.p. functions on $G$ is denoted by $WAP(G)$. Note that $WAP(G)$ is a C*-subalgebra $C(G)$. If $G$ is compact, then $C(G) = C_0(G) = AP(G) = WAP(G)$. If $G$ is noncompact then $AP(G) \oplus C_0(G) \subset WAP(G)$. In this note, we usually will only be interested in noncompact groups $G$. The algebra $WAP(G)$ was first introduced and studied by Eberlein [11] when $G$ is abelian; Burckel [2] is a convenient reference for many of the earlier results on weakly almost periodic functions.

In the mid 1980's, I noticed that the double orbits of weakly almost periodic functions may not be relatively weakly compact. I never published my findings on double orbits of w.a.p. functions but did share them to a few researchers whose research interests are close to mine. One of my initial examples is the following: Let $M(2)$ be the two-dimensional motion group. If $f \in C_0(M(2)), f \neq 0$, then $O(f)$ is not relatively weakly compact. This example was given as Exercise 2.24 on p. 149 of the monograph [1]; see also the comments on p. 218 of [1]. It was also mentioned on p. 345 of Lau and Ülger [22].

Independently, G. Hansel and J.P. Troallic provided a more systematic study of the double orbits of w.a.p. functions in a sequence of three papers in the early 1990's; see [17], [18] and [19]. In this note we will adapt their terminologies:

**Definition 1.1.** For a locally compact group $G$, let
$$WS(G) = [f \in C(G) : O(f) \text{ is relatively weakly compact in } C(G)\}.$$

Functions in $WS(G)$ are said to be strictly weakly almost periodic. Note that

$$AP(G) \subset WS(G) \subset WAP(G).$$

As in [17], we will call $G$ a $WS$-group, or $G \in [WS]$, if $WS(G) = WAP(G)$.

Clearly, abelian groups and compact groups are $WS$-groups. Here is a main result of [17]:

**Theorem 1.2.** ([17], Theorem 4.3) The following 2 conditions are equivalent: (1) $C_0(G) \subset WS(G)$. (2) The left and right uniform structures on $G$ are equal.

Recall that the left and right uniform structures of a locally compact group $G$ are equal if and only if $G$ is a $SIN$-group, i.e., the collection of neighborhoods of the identity $e$ of $G$, invariant under the inner automorphisms of $G$, forms a neighborhood basis at the identity of $G$; see Hewitt and Ross [20], p. 21. So being a $SIN$-group is a necessary condition for a group to be a $WS$ group. But it is not a sufficient condition. We will give, in this note, examples of discrete groups which are not $WS$-groups.

The main result of [19] identifies the currently known $WS$-groups:

**Theorem 1.3.** ([19], Theorem 4.2) If $G$ is a locally compact $Moore$-group then $G \in [WS]$.

They asked whether $G \in [WS]$ would imply that G is a $Moore$-group. This problem appears to be still unsolved. Recall that a locally compact group $G$ is call a $Moore$-group, if all irreducible continuous unitary representations of $G$ are finite-dimensional; see Moore [25]. We asked in [7] whether a discrete $WS$-group must be a finite extension of an abelian group. A classical result of Thoma [31], states that a discrete group is a $Moore$-group, if and only if it is abelian by finite. Therefore, as was also pointed out in [19], our question is the restriction of their question to discrete groups. We will provide a positive answer to this question for a smaller class of discrete groups in Section 4: a discrete $FC$-group is a $WS$-group if and only it is a finite extension of its center. Recall that a group $G$ is an $FC$-group if each conjugacy class of $G$ is finite and there are known examples of discrete $FC$-groups which are not finite extensions of abelian groups.

In Section 3, we will study noncompact locally compact groups $G$ which satisfy the condition that $WS(G) = AP(G)$; i.e., for $f \in C(G)$, whenever $O(f)$ is relatively weakly compact then it is relatively norm compact. We will show that $WS(G) \neq AP(G)$ if $G$ is either a noncompact $IN$-group or a noncompact nilpotent group. On the other hand, when $G$ is minimally w.a.p. then $WS(G) = AP(G)$. Recall that $G$ is called a minimally w.a.p. group if $WAP(G) = AP(G) \oplus C_0(G)$; see Chou [5]. We showed in 1975 [4, Theorem 4.8] that $M(2)$ is minimally w.a.p. We also showed in 1980 [5, Theorem 3.1] that if $G$ is a connected solvable minimally w.a.p. group and $K(G)$ is the largest compact normal subgroup of $G$ then $G/K(G)$ is topologically isomorphic to $M(2)$.

To study strictly w.a.p functions, it is convenient to introduce the following.

**Definition 1.4.** Let $UCS(G) = \{f \in C(G): O(f)$ is equicontinuous$\}$.

Note that $WS(G) \subset UCS(G)$. We will show that $UCS(SL(2,\mathbb{R}))$ only contains constant functions.

## 2. Preliminaries and $IN$-groups

Let $G$ be a locally compact group. $f \in C(G)$ is said to be left uniformly continuous, if given $\varepsilon > 0$, there exists a neighborhood $U$ of the identity $e$ of $G$ such that $|f(s) - f(t)| < \varepsilon$, whenever $st^{-1} \in U$; i.e., $f$ is uniformly continuous with respect to the right uniform structure of $G$. Note that $f$ is left uniformly continuous, if $x \to {}_xf$ is continuous from $G$ to $C(G)$. Let $LUC(G)$ be the $C^*$-algebra of all bounded left uniformly continuous functions on $G$. Similarly, we can define $RUC(G)$, the algebra of all bounded right uniformly continuous functions on $G$ and $UC(G) = LUC(G) \cap RUC(G)$, the algebra of (two-sided) uniformly continuous functions on $G$; see Hewitt and Ross [19]. Clearly, $UCS(G) \subset UC(G)$ and it is known that $WAP(G) \subset UC(G)$; see [2].

**Lemma 2.1.** Let $G$ be a locally compact group.
(1) If $f \in LUC(G)$ or $RUC(G)$ and $\{{}_xf_{x^{-1}}: x \in G\}$ is equicontinuous at $e$ then $f \in UCS(G)$.
(2) $G$ is a $SIN$-group if and only $UCS(G) = UC(G)$.
(3) $WS(G) \subset UCS(G)$.

**Proof.** (1) Note that, for $f \in C(G)$, $f(xuy) - f(xy) = f(xux^{-1}xy) - f(xy)$ for $x, y, u \in G$.

(2) This is part of Lemma 4.1 of Hansel and Troallic [17].

(3) In the proof of (1) $\Rightarrow$ (2) of Theorem 4.3 of [17], Hansel and Troallic showed that if $f \in WS(G)$, then $f \in UCS(G)$, using Robert Ellis' joint continuity theorem; see [1]. ∎

(However, Hansel and Troallic only stated (3) of the above lemma for functions in $C_0(G)$, assuming $C_0(G) \subset WS(G)$.)

Lau and Ülger gave a different proof of the fact that if $G$ is a $SIN$-group, then $C_0(G) \subset WS(G)$, using the fact that the von Neumann algebra of a $SIN$-group is finite; see [22, Proposition 7.16].

If $N$ is a closed normal subgroup of a locally compact group $G$, we will denote the coset $xN$ by $\dot{x}$. If $N$ is, in addition, compact, for $f \in UC(G)$, let $f^N(\dot{x}) = \int_N f(xt)dt$. Here the integral is with respect to the normalized Harr measure on $N$.

**Lemma 2.2.** If $N$ is a compact normal subgroup of a locally compact Group $G$ and If $f \in UCS(G)$, then $f^N \in UCS(G/N)$.

**Proof.** Note that $f^N(\dot{x}\dot{u}\dot{y}) = \int_N f(xuyt)dt$. ∎

A locally compact group $G$ is called an $IN$-group, if it has a compact invariant neighborhood of the identity $e$. It is known that if $G$ is an $IN$-group then the intersection of all closed invariant neighborhoods of $e$ is a compact normal subgroup $K = K_G$ of $G$ and the quotient group $G/K$ is a $SIN$-group; see Iwasawa [21].

**Lemma 2.3.** Assume that $G$ is an $IN$-group; let $K = K_G$ be the compact normal subgroup of $G$ defined above and let $\theta$ be the natural homomorphism of $G$ onto $G/K$. Then $UCS(G) = \{h \circ \theta : h \in UCS(G/K)\}$.

**Proof.** Let $f \in UCS(G)$. We claim that $f$ is constant on the cosets of $K$. Indeed, for $\varepsilon > 0$, let $W_\varepsilon = \{x \in G : |f(x) - f(e)| \leq \varepsilon\}$. Note that $W_\varepsilon$ is a closed neighborhood of $e$. Since $f \in UCS(G)$, the set of functions $\{_xf_{x^{-1}} : x \in G\}$ is equicontinuous at $e$. So, there is a closed neighborhood $V_\varepsilon$ of $e$ such that if $u \in V_\varepsilon$ and $x \in G$, then $|f(xux^{-1}) - f(e)| \leq \varepsilon$. Therefore, if $u \in V_\varepsilon$ then $xux^{-1} \in W_\varepsilon$, and hence $V_\varepsilon \subset \cap \{x^{-1}W_\varepsilon x : x \in G\} = U_\varepsilon$. So, $U_\varepsilon$ is a closed invariant neighborhood of $e$. By the definition of $K$, $K \subset U_\varepsilon \subset W_\varepsilon$. So, if $t \in K$, then $|f(t) - f(e)| \leq \varepsilon$. Since $\varepsilon > 0$ is arbitrary, $f(t) = f(e)$, if $t \in K$. Let $xK$ be a coset of $K$ in $G$. Note that since $_xf \in UCS(G)$, $f(xt) = {_xf}(e) = f(x)$, for all $t \in K$. Our claim is proved.

Let now $f \in UCS(G)$. Then, by Lemma 2.2, $f^K \in UCS(G/K)$. Since $f$ is constant on cosets of $K$, $(f^K) \circ \theta = f$. ∎

For a general locally compact group $G$, $WAP(G)$ has a unique invariant mean, denoted by $m$ or $m_G$; see Ryll-Nardzewski [32]. Let $WAP_0(G) = \{f \in WAP(G) : m(|f|) = 0\}$ then $WAP(G) = AP(G) \oplus WAP_0(G)$; see [1]. When $G$ is noncompact, $C_0(G) \subset WAP_0(G)$. If $C_0(G) = WAP_0(G)$ then $G$ is called a minimally w.a.p. group; see [5]. For example, $M(n)$, the $n$ dimensional motion group, $n \geq 2$, and noncompact simple analytic groups with finite centers are minimally w.a.p.; see Chou [5], Veech [35].

Let $B(G)$ be the Fourier-Stieltjes algebra of $G$ and $B(G)^-$ be its uniform closure in $C(G)$. When $G$ is abelian, $B(G)$ is the algebra of all Fourier-Stieltjes transforms of bounded regular Borel measures on the dual group $\hat{G}$ of $G$. For a general locally compact group, $B(G)$ is the algebra of coefficient functions of continuous unitary representations of $G$ and it was first defined and studied by Eymard [13]. It is known that $B(G)^- \subset WAP(G)$ and if $G$ is compact then $B(G)^- = C(G)$. Eberlein raised the question whether, for a noncompact abelian group $G$, $B(G)^- = WAP(G)$. The answer turned out to be negative for all noncompact abelian groups; see Rudin [29] and Ramirez [28]. In [6] we called $G$ an Eberlein group, if $B(G)^- = WAP(G)$. So, noncompact abelian groups are not Eberlein groups. We extended their results to many nonabelian groups in [6]: if $G$ is a noncompact $IN$-group or a noncompact nilpotent group then $G$ is not an Eberlein group; in fact, the quotient Banach space $WAP(G)/B(G)^-$ contains an isometric copy of $\ell^\infty$. More recently, Filali and Galindo [14] were able to show that the quotient space for these two classes of locally compact groups contains an isometric copy of $\ell^\infty(\kappa)$ where $\kappa$ is the minimal number of compact sets required to cover $G$; see Theorems 5.6 and 5.7 of [14].

If $N$ is a closed (not necessarily compact) normal subgroup of a locally compact group $G$, then for $f \in WAP(G), x \in G$, let $f^N(\dot{x}) = m_N(f^x)$ where $f^x \in C(N)$ is defined by $f^x(t) = f(xt), t \in N$. Since $m_N$ is translation invariant on $WAP(N)$, $f^N$ is well defined. In fact, $f^N \in WAP(G/N)$; see Lemma 2.3 of Chou [5].

**Lemma 2.4.** Let $\tau$ be a continuous automorphism of a locally compact group $G$. If $f \in WAP(G)$, then $f \circ \tau \in WAP(G)$ and $m(f \circ \tau) = m(f)$ where $m$ is the unique invariant mean on $WAP(G)$.

Proof. The fact that $f \circ \tau$ is w.a.p. is a direct consequence of Grothendieck's weak compactness criterion. For $f \in WAP(G)$, let $m'(f) = m(f \circ \tau)$. Clearly $m'$ is a mean on $WAP(G)$. Note that $(_x f) \circ \tau = {}_{\tau^{-1}(x)}(f \circ \tau)$. So, $m'(_x f) = m(f \circ \tau) = m'(f)$, i.e., $m'$ is a left invariant mean on $WAP(G)$. By the uniqueness of invariant mean on $WAP(G)$, $m' = m$. ∎

**Lemma 2.5.** Assume $N$ is a closed normal subgroup of a locally compact group $G$. If $f \in WS(G)$ then $f^N \in WS(G/N)$.

**Proof.** Note that for $a, b \in G$,

(2.1) $\qquad (_a f_b)^N = {}_{\dot{a}}(f^N)_{\dot{b}}.$

Indeed, for $x \in G$ and $t \in N$,
$$(_a f_b)^x(t) = (_a f_b)(xt) = f(axtb) = f(axbb^{-1}tb) = f^{axb}(b^{-1}tb).$$

For a fixed $b \in G$, $\tau: t \to b^{-1}tb$ is a continuous automorphism of $N$. By Lemma 2.4, $m_N(f^{axb} \circ \tau) = m_N(f^{axb})$.
So, $(_a f_b)^N(\dot{x}) = m_N(((_a f_b)^x) = m_N(f^{axb} \circ \tau) = m_N(f^{axb}) = f^N(\dot{a}\dot{x}\dot{b})$ and we have proved (2.1). To complete the proof of this lemma, just follow the steps of Lemma 2.3 of [5]. ∎

**Theorem 2.6.** Let $G$ be a locally compact $IN$-group. Let $K_G$ and $\theta$ be defined as in Lemma 2.3. Then $WS(G) = \{h \circ \theta : h \in WS(G/K_G)\}$.

**Proof.** Let $f \in WS(G)$. Then, by Lemma 2.1(3), $f \in UCS(G)$. Lemma 2.3 implies that $f$ is constant on the cosets of $K = K_G$ in $G$ and $f = f^K \circ \theta$. By Lemma 2.5, $f^K \in WS(G/K)$. ∎

From now on if $N$ is a closed normal subgroup of a locally compact group and $\theta$ is the natural homomorphism of $G$ onto $G/N$, we often identify $g \in C(G/N)$ with $g \circ \theta \in C(G)$. For example, Theorem 2.6 states that $WS(G) = WS(G/K_G)$.

It is known that if $G$ is a connected $IN$-group, then $G/K_G$ is a connected $SIN$-group and hence is a direct product of a vector group $R^n$ and a compact group; see Grosser and Moskowitz [15]. Therefore, $G/K_G$ is a $WS$-group, and hence, by Theorem 2.6 and Lemma 2.3, we have the following:

**Corollary 2.7.** Let $G$ be a connected $IN$-group. Then $WS(G) = WAP(G/K_G)$ and $UCS(G) = UC(G/K_G)$.

**Example 2.8.** Let $G = \left\{ \begin{bmatrix} 1 & x & e^{it} \\ 0 & 1 & y \\ 0 & 0 & 1 \end{bmatrix} : x, y, t \in \mathbb{R} \right\}$ be the reduced Heisenberg group. Then $G$ is a rank 2 nilpotent group. We will write elements of $G$ as $(x, y, e^{it})$. Then the center of $G$ is $K =$ $\{(0,0,e^{it}), t \in \mathbb{R}\} \simeq \mathbb{T}$ and $G$ is an $IN$-group and $K$ is the intersection of all closed invariant neighborhoods of $e$. By Corollary 2.7, $UCS(G) = UCS(G/K) = UC(\mathbb{R}^2)$ and $WS(G) = WS(G/K) = WAP(\mathbb{R}^2)$. Note that here $\mathbb{R}^2$ is considered as a quotient of $G$ and the subset $\{(x,y,1): x, y \in \mathbb{R}\}$ is not a subgroup of $G$. It is known that $WAP(G)|K = WAP(K) = C(K)$; see the discussions on p. 92 of Cowling and Rodway [8]. But $WS(G)|K \neq WS(K)$. In fact, $WS(G)|K$ only contains constant functions.

The following proposition is similar to Proposition 7.18 of Lau and Ülger [21].

**Proposition 2.9.** The following three conditions on a locally compact group $G$ are equivalent: (1) $G$ is an $IN$-group; (2) $C_0(G) \cap WS(G) \neq \{0\}$; (3) $C_0(G) \cap UCS(G) \neq \{0\}$.

**Proof.** (1) $\Rightarrow$ (2). Assume that $G$ is an $IN$-group. Let $K$ be the intersection of all compact invariant neighborhoods of $e$ in $G$. Then $G/K$ is a $SIN$-group, and by Theorem 1.2, $C_0(G/K) \subset WS(G/K)$. Choose any $h \in C_0(G/K), h \neq 0$. Then, by Theorem 2.6, $h \circ \theta \in C_0(G) \cap WS(G)$.
(2) $\Rightarrow$ (3) is obvious, since by Lemma 2.1(3), $WS(G) \subset UCS(G)$.
(3) $\Rightarrow$ (1). Choose $f \in C_0(G) \cap UCS(G)$ and $f \neq 0$. We may assume that $f(e) = 1$. Let $W = \{x \in G: |f(x) - 1| \leq \frac{1}{2}\}$. Then $W$ is a compact neighborhood of $e$. Since $\{_x f_{x^{-1}}: x \in G\}$ is equicontinuous at $e$, $U = \cap \{xWx^{-1}: x \in G\}$ is a compact invariant neighborhood of $e$ and hence $G$ is an $IN$-group. ∎

In the sequel, $G$ will denote a locally compact group and a subset $U$ of $G$ is said to be invariant will mean that $U$ is invariant under the inner automorphisms of $G$.

### 3. $WS(G)$ and $AP(G)$

**Lemma 3.1.** (1) Assume $N$ is a closed normal subgroup of $G$. If $WS(G) = AP(G)$, then $WS(G/N) = AP(G/N)$. (2) If $K$ is a compact normal subgroup of $G$, and $WS(G/K) = AP(G/K)$, then $WS(G) = AP(G)$.

**Proof.** (1) Assume $WS(G/N) \supsetneq AP(G/N)$. Let $g \in WS(G/N) \setminus AP(G/N)$. Consider the Eberlein decomposition of $g$ in $WAP(G/N)$: $g = g_1 + g_2$, $g_1 \in WAP_0(G/N)$, $g_2 \in AP(G/N) \subset WS(G/N)$. Therefore, $g_1 \in WS(G/N) \cap WAP_0(G/N)$ and $g_1 \neq 0$. Then $g_1 \circ \theta \in WS(G)$, but $g_1 \circ \theta \notin AP(G)$.

(2) Assume $AP(G) \subsetneq WS(G)$. As in (1), by considering Eberlein decomposition, we may assume that there exists $f \in WS(G) \cap WAP_0(G)$. We may assume $f \geq 0$ and $f \neq 0$. By Lemma 2.5, $f^K \in WS(G/K)$, by Lemma 2.3 of [5], $f^K \in WAP_0(G/K)$. Since $K$ is compact, $f \geq 0, f \neq 0$ and $f^K(\dot{x}) = \int_K f(xt)dt$, $f^K \neq 0$. Hence, $f^K \notin AP(G/K)$. ∎

**Theorem 3.2.** (1) If $G$ is a noncompact locally compact $IN$-group, then $WS(G) \supsetneq AP(G)$. (2) If $G$ is a noncompact nilpotent groups then $WS(G) \supsetneq AP(G)$.

**Proof.** (1) Assume $K$ is the intersection of all closed invariant neighborhoods of the noncompact $IN$-group $G$. Then $C_0(G/K) \subset WS(G/K) \cap C_0(G) = WS(G) \cap C_0(G)$. Then, clearly, $WS(G) \supsetneq AP(G)$.

(2) Assume $G$ is a noncompact nilpotent group. Consider the upper central sequence of $G$:
$$G = G_0 \supset G_1 \supset \cdots G_{n-1} \supset G_n = (e).$$
Where each $G_i$ is a closed normal subgroup of $G$ and $G_{i-1}/G_i$ is the center of $G/G_i$. We will call $n$ (the nilpotent) rank of $G$.

If $n = 1$, then $G$ is abelian and hence $WS(G) = WAP(G) \supsetneq AP(G)$. Assume $n > 1$ and (2) holds for all noncompact nilpotent groups of rank $\leq n - 1$. Now $G_{n-1} = Z(G)$, the center of $G$. If $G/Z(G)$ is compact, then $G$ is called a $Z$-group and hence $WS(G) = WAP(G) \supsetneq AP(G)$; see Corollary 5.3 of [18]. If $G/G_{n-1}$ is noncompact then by inductive assumption, $WS(G/G_{n-1}) \supsetneq AP(G/G_{n-1})$. Then, Lemma 3.1(1) implies that $WS(G) \supsetneq AP(G)$. ∎

Note that Theorem 3.2 (2) cannot be extended to solvable locally compact groups. For example, consider the solvable group $G = M(2)$. As mentioned in Section 1, $C_0(G) \cap WS(G) = (0)$. Since $G$ is minimally w.a.p., $WAP(G) = AP(G) \oplus C_0(G)$ and therefore, $WS(G) = AP(G)$. We will now show that all minimally w.a.p. groups have this property.

**Theorem 3.3**. Let $G$ be a noncompact minimally weakly almost periodic group. Then $WS(G) = AP(G)$.

**Proof.** By the von Neumann approximation theorem for almost periodic functions, $AP(G) \subset B(G)^-$. It is also well-known that the Fourier algebra $A(G)$ of $G$, which is a subalgebra of $B(G)$, is uniformly dense in $C_0(G)$; see Eymard [13]. Since $G$ is minimally w.a.p., $WAP(G) = AP(G) \oplus C_0(G) \subset B(G)^-$. Therefore, $G$ is an Eberlein group; by Theorem 4.5 of [6], $G$ is not an $IN$-group. By Proposition 2.9, $C_0(G) \cap WS(G) = (0)$. Hence, $WS(G) = AP(G)$. ∎

We showed in [5] that $SL(2, \mathbb{R})$ and $M(n) = \mathbb{R}^n \rtimes SO(n)$, the $n$-dimensional motion group are minimally w.a.p. More generally, Veech [33] proved that simple analytic groups with finite centers are minimally w.a.p. By applying the theory of totally minimal topological groups, Mayer [23], showed that $G = \mathbb{R}^n \rtimes K$ is minimally w.a.p., if $K$ is a compact group acting on $\mathbb{R}^n$ irreducibly.

Recall that a locally compact group $G$ with topology $\tau$ is said to be minimal if $\tau$ contains no strictly coarser Hausdorff topologies (not necessarily locally compact) which make $G$ a topological group. $G$ is totally minimal if the quotient group $G/N$ is minimal for every closed normal subgroup $N$ of $G$; see Dikranjan, Prodanov, and Stoyanov [10]. For a more recent survey on totally minimal groups, see Dikranjan and Megrelishvili [9].

The following theorem of Mayer [24], gives necessary and sufficient conditions for a connected locally compact group to be an Eberlein group which summarizes his main findings in [23]:

**Theorem 3.4** ([24], Theorem 5). The following conditions on a connected locally compact group $G$ are equivalent:

(1) $G$ is totally minimal.
(2) There exists a compact normal subgroup $K$ such that the quotient group $G/K = N \rtimes H$, a semidirect product of $N$ and $H$ where $N$ is a simply connected nilpotent analytic group and $H$ is a linear reductive group and $H$ acts on $N$ without nontrivial fixed points.
(3) $G$ is an Eberlein group.
(4) $WAP(G)$ is the $C^*$-algebra generated by functions in $C_0(G/N)$ where $N$ ranges over all closed normal subgroups of $G$.

By the above theorem, if $G$ is a connected Eberlein group, then $G$ is totally minimal and, by (4), $WAP(G)$ is generated as a $C^*$-algebra by (the pull backs of) functions in $C_0(G/N)$ when $N$ ranges over all closed normal subgroups of $G$. There are two possible cases:

(i) $G/N$ is compact. Then $AP(G/N) = C_0(G/N) = C(G/N) \subset AP(G)$.

(ii) $G/N$ is noncompact. Then $G/N$ is a noncompact totally minimal group. By the above theorem, $G/N$ is a noncompact Eberlein group. Therefore, by Theorem 4.5 of [6], $G$ is not an $IN$-group. By Proposition 2.9, $C_0(G/N) \cap WS(G/N) = (0)$. Note also $C_0(G/N) \subset WAP_0(G)$.

Because of the Eberlein decomposition $WAP(G) = AP(G) \oplus WAP_0(G)$ and (4) of Theorem 3.4, we wonder if $G$ is a connected Eberlein group then $WAP_0(G) \cap WS(G) = (0)$, i.e., $WS(G) = AP(G)$. Using Theorem 3.4 (2), Mayer provided in [24] the following examples of Eberlein groups: (1) the Lorentz groups $\mathbb{R}^n \rtimes SL(n, \mathbb{R})$, (2) the Euclidean motion groups, and (3) $(\mathbb{R}^2 \oplus \mathbb{R}^3) \rtimes (SL(2, \mathbb{R}) \times SL(3, \mathbb{R}))$. Note that, by [5, Lemma 2.2(b)], for $n \geq 2$, the $n$-dimentional Lorentz group is not minimally w.a.p., since it is a semidirect product of two noncompact groups.

W. Veech proved in [35] that if $G$ is a noncompact simple analytic group with finite center, then $WAP(G) = \mathbb{C} \oplus C_0(G)$. In particular $G$ is minimally w.a.p. Let $G$ be a semisimple analytic group with finite center. Then there is a finite extension $G_0$ of $G$ such that $G_0$ is a direct product

(3.1) $$G_0 = G_1 \times G_2 \times \ldots \times G_n$$

where each $G_i$ is a simple analytic group with finite center.

**Theorem 3.5.** Let $G$ be a semisimple analytic group with finite center. Then $WS(G) = AP(G)$.

**Proof.** Assume $G$ is a semisimple analytic group with finite center, as was just described above. To show that $WS(G) = AP(G)$, by Lemma 3.1, we may assume $G = G_0$ in (3.1) and by factoring out the compact factors, we may further assume that all the $G_i$'s are noncompact. Let $f \in WS(G)$ and $x \in G$. We will write $f(x) = f(x_1, x_2, \ldots, x_n)$ where $x_i \in G_i, i = 1, 2, \ldots, n$. Fix $x_2, \ldots, x_n$ and let $g(t) = f(t, x_2, \ldots, x_n), t \in G_1$. Then $g \in WS(G_1)$. Since $G_1$ is minimally w.a.p., by Theorem 3.3, $g \in AP(G_1) = \mathbb{C}$, the constant functions on $G_1$. By similar considering for the other variables, we conclude that any $f \in WS(G)$ is a constant function with respect to the $i$th variable when the remaining variables are fixed. This, of course, implies that $f$ is a constant function on $G$. ∎

## 4. $WS$-groups

The main theorem of [19] states that locally compact $Moore$-groups are $WS$-groups. It is not known whether the converse holds. The known structure theorems of almost connected groups imply the following.

**Proposition 4.1.** Let $G$ be an almost connected locally compact group. Then the following conditions on $G$ are equivalent: (1) $G$ is $WS$-group; (2) $G$ is a $SIN$-group; (3) $G$ is a $Moore$-group; (4) $G = V \rtimes_\varphi N$, a semidirect of a vector group $V$ and a compact group $N$ where $\varphi(N)$ is finite.

**Proof.** (1) $\Rightarrow$ (2). If $G$ is a $WS$-group, then $C_0(G) \subset WS(G)$, and, by Theorem 1.2, $G$ is a $SIN$-group. (2) $\Rightarrow$ (1). See, for example, Corollary 5.3 of [17].

The structure theorems of locally compact groups imply that (2), (3) and (4) are equivalent; see T. W. Palmer's survey article [29]. ∎

For the remainder of this section, we will only consider discrete groups. Let $G$ be a discrete group. A subset $S$ of $G$ is called a $t$-set if the sets $xS \cap S$ and $Sx \cap S$ are finite whenever $x \in G$ and $x \neq e$. We will need the following:

**Lemma 4.2.** (1) Every infinite subset of a group $G$ contains an infinite $t$-set.
(2) Assume $S$ is a $t$-set in $G$. Then if $f \in \ell^\infty(G)$ and $f(x) = 0$ whenever $x \notin S$ then $f \in WAP(G)$.

Both assertions of the above lemma are known, for example, see Proposition 4.1 of [3] for (1) and Lemma 3.2 of [6] for (2).

For $x \in G$, the conjugacy class of $G$ containing $x$ is denoted by $cl(x)$ or $cl_G(x)$. $G$ is called an $FC$-group if each conjugacy class of $G$ is finite. The monograph of Tomkinson [34], is a convenient place to look for basic facts of $FC$-groups. We will adapt the terminology of [34] to call a subgroup of $G$ an $N$-subgroup if it is generated by elements $a_i, b_i, i \in \mathbb{N}$, and the generators satisfy the following conditions:

(4.1) $\quad [a_i, a_j] = [b_i, b_j] = [a_i, b_j] = 1, \text{if } i \neq j; [a_i, b_i] = c_i \neq e.$

As usual, for $a, b \in G$, $[a, b] = aba^{-1}b^{-1}$. To prove the main result of this section, we need the follow theorem of B. H. Neumann [27]:

**Theorem 4.3.** (Neuman [27]) Let $G$ be an $FC$-group. If $Z(G)$, the center of $G$, is of infinite index in $G$, then $G$ contains an $N$-subgroup.

**Theorem 4.4.** Let $G$ be an $FC$-group. Then $G$ is a $WS$-group, i.e., $WS(G) = WAP(G)$, if and only $G$ is a finite extension of its center.

**Proof.** It is easy to prove directly that finite extensions of abelian groups are $WS$-groups. It is, of course, also a consequence of Theorem 4.2 of [19].

Assume that $G$ is an $FC$-group and $Z(G)$ is of infinite index in $G$. Then by Theorem 4.3, $G$ contains sequences of elements $a_i, b_i$ satisfying (4.1). Applying Lemma 4.1(1), by taking corresponding subsequences of $(a_i), (b_i)$, we may assume $B = \{b_i : i = 1,2, ...\}$ is a $t$-set, and by Lemma 4.1 (2), $\chi_B$, the characteristic function of $B$ in $G$, belongs to $WAP(G)$.

Let $x_n = a_1 a_2 ... a_n$. Then

(4.2) $$x_n b_m x_n^{-1} = \begin{cases} c_m b_m, & n \geq m \\ b_m, & n < m. \end{cases}$$

Note that $c_m b_m \neq b_k$, if $k \in \mathbb{N}$. Indeed, $c_m b_m \neq b_m$, since $c_m \neq e$. If $k \neq m$ and $c_m b_m = b_k$, then $a_m b_m a_m^{-1} b_m^{-1} b_m = b_k$, and hence, $b_m = a_m^{-1} b_k a_m = b_k$, a contradiction. Therefore,

$$\lim_n \lim_m \chi_B(x_n b_m x_n^{-1}) = 1$$

but

$$\lim_m \lim_n \chi_B(x_n b_m x_n^{-1}) = 0.$$

By Grothendieck's criterion, $\chi_B \notin WS(G)$. ∎

The above proof shows that if $D$ is an infinite subset of $B$ then $\chi_D \in WAP(G) \setminus WS(G)$.

Since a finite extension of an abelian group may not be an $FC$-group, it is easy to find $WS$-groups which are not $FC$-groups. For example, the infinite dihedral group $D_\infty \in [WS]$, but it is not an $FC$-group. We would like to include two examples of discrete non-$WS$-groups here.

**Examples 4.5.**

(1) (J. Erdös [12, p. 58]) Let $p$ be a fixed prime number and let $G_1$ be the group with infinite many generators: $b, a_i, i \in \mathbb{N}$, and with defining relations: for $i, k \in \mathbb{N}$,
$$b^p = a_i^p = e$$
$$a_i b = b a_i$$
$$a_{i+k} a_i = b a_i a_{i+k}.$$
Then $G_1'$, the derived group of $G_1$, is the cyclic group $\langle b \rangle$ of order $p$. Hence $|cl(x)| \leq p$ for $x \in G_1$ and hence $G_1$ is an $FC$-group. Note that $Z(G_1) = \langle b \rangle$ is finite and hence is of infinite index in $G_1$. By the above theorem, $G_1$ is not a $WS$-group. (Note that $G_1$ is infinitely generated. It is known and is easy to show that the center of a finitely generated $FC$-group $G$ is always of finite index in $G$.)

(2) Let $G_2 = N \rtimes_\eta H$ where $N = \mathbb{Z}\left[\frac{1}{2}\right] = \{\frac{m}{2^k}: m, k \in \mathbb{Z}\}$ and $H = \mathbb{Z}$ and the action of $H$ on $N$ is given by $\eta(n)\left(\frac{m}{2^k}\right) = 2^n \frac{m}{2^k}$. Consider

(4.3) $\qquad (0, -n)(2^m, 0)(0, n) = (2^{m-n}, 0).$

Note that $N_1 = \{(m, 0): m \in \mathbb{Z}\}$ is a subgroup of $N \times (0)$, $N_1 \simeq \mathbb{Z}$ and $S = \{(2^m, 0): m \in \mathbb{N}\}$ is a $t$-subset of $N_1$, and hence is a $t$-subset of $G_2$. By Lemma 4.1 (2) $\chi_S \in WAP(G_2)$. By (4.3),
$$\lim_n \lim_m \chi_S[(0, -n)(2^m, 0)(0, n)] = 1$$
but
$$\lim_m \lim_n \chi_S[(0, -n)(2^m, 0)(0, n)] = 0.$$
By Grothendieck's criterion, $\chi_S \notin WS(G)$. So $G_2$ is not a $WS$-group. (Note that $G_2$ is the Baumslag–Solitar Group $B(1,2)$. It is a finite presented solvable group.)

## 5. The $ax + b$ group and $SL(2, \mathbb{R})$

We showed in Lemma 2.3 and Theorem 2.6 that if $G$ is an $IN$-group and $K$ is the intersection of all closed invariant neighborhoods of $e \in G$, then $UCS(G) = UCS(G/K)$ and $WS(G) = WS(G/K)$. We will now describe $WS$ functions on the $ax + b$ group which is not an $IN$-group.

Let $G_3 = \mathbb{R} \rtimes \mathbb{R}^+$ be the $ax + b$ group. Recall the multiplication on $G_3$ is
$$(b, a)(b', a') = (b + ab', aa').$$
Note that $(b, a)^{-1} = (-\frac{b}{a}, \frac{1}{a})$ and

(5.1) $\qquad (b, a)(b', a')(b, a)^{-1} = (ab' + b(1 - a'), a').$

**Lemma 5.1.** If $V$ is an invariant neighborhood of $e = (0,1) \in G_3$, then there exists $0 < \varepsilon < 1$, such that $V_\varepsilon = \{(d, c): d \in \mathbb{R}, |c - 1| < \varepsilon\} \subset V$.

**Proof.** Because $V$ is a neighborhood of $e$, there exists $0 < \varepsilon < 1$ such that
$$U = \{(b', a'): |b'| < \varepsilon, |a' - 1| < \varepsilon\} \subset V.$$
Since $V$ is invariant,
$$W = \bigcup_{(b,a) \in G_3} (b, a)U(b, a)^{-1} \subset V$$
Claim that $W = V_\varepsilon$. Indeed, let $(b', a') \in U$ and $x \in \mathbb{R}$.

If $a' \neq 1$, by (5.1),
$$\left(\frac{x}{1 - a'}, 1\right)(0, a')\left(\frac{x}{1 - a'}, 1\right)^{-1} = (x, a').$$
If $a' = 1, x \neq 0$, choose any $b' \in \mathbb{R}, 0 < |b'| < \varepsilon$. Then

$$\left(0, \frac{x}{b'}\right)(b', 1)\left(0, \frac{x}{b'}\right)^{-1} = (x, 1).$$

The claim is proved. ∎

**Theorem 5.2.** Let $G_3$ be the $ax + b$ group defined above.

(1) If $f \in UCS(G_3)$, then it is constant on each of the coset of $N = \mathbb{R} \times 1$ in $G_3$.
(2) $WS(G_3) = WAP(G_3/N) = WAP(\mathbb{R}^+)$.

**Proof.** (1) Let $f \in UCS(G_3)$. Then $\{_x f_{x^{-1}} : x \in G\}$ is equicontinuous at $e$ by following the proof of Lemma 2.3, one sees that, for each $k \in \mathbb{N}$, there is an invariant neighborhood $W_k$ of e in $G_3$ such that if $z \in W_k$ then $|f(z) - f(e)| < \frac{1}{k}$. By Lemma 5.1, there exists $\varepsilon_k$, $0 < e_k < 1$, such that $V_{\varepsilon_k} \subset W_k$. Note that if $z \in \cap_k V_{\varepsilon_k}$, then $f(z) = f(e)$. Clearly, we may choose $\varepsilon_k$ to be a decreasing sequence and $\lim_k \varepsilon_k = 0$. Then $\cap_k V_{\varepsilon_k} = N$. So $f(z) = f(e)$, if $z \in N$. By replacing $f$ by $_z f$, for $z \in G_3$, we conclude that $f$ is constant on each of the coset of $N$ in $G_3$.

(2) Let $f \in WS(G_3)$. By Lemma 2.1(3), $f \in UCS(G_3)$ and hence, by (1), $f$ is constant on each of the coset of $N$. Define $f^N$ as in Lemma 2.5. Then $f^N \in WS(G/N)$, Since $f$ is constant on cosets of $N$ in $G_3$, $f = f^N \circ \theta$, where $\theta$ is the natural homomorphism of $G_3$ onto $G_3/N$. ∎

**Remarks.** (1) Let $G = M(2) = \mathbb{C} \rtimes \mathbb{T}$ be the 2-dimensional motion group. Recall the multiplication on $M(2)$ is given by $(z, w)(z'w') = (z + wz', ww')$. Let $N = \mathbb{C} \times 1$. Then using the same proof as that of Theorem 5.2, one can conclude that
$$WS(M(2)) = WAP(M(2)/N) = WAP((0) \times \mathbb{T}) = AP(\mathbb{T}) = AP(M(2)).$$
This also follows from the fact that $M(2)$ is minimally w.a.p. But the approach outlined here is more direct and simpler.

(2) P. Milnes studied w.a.p. functions on $M(2)$, the $ax + b$ group and the reduced Heisenberg group in [24].

**Theorem 5.3.** $UCS(SL(2, \mathbb{R}))$ only contains constant functions.

**Proof.** For convenience, we will denote the group $SL(2, \mathbb{R})$ by $G$. Consider the Iwasawa decomposition of $G$: $G = KAN$, where
$$K = \left\{\gamma(\theta) = \begin{bmatrix} \cos\theta & \sin\theta \\ -\sin\theta & \cos\theta \end{bmatrix} : 0 \leq \theta \leq 2\pi\right\},$$
$$A = \left\{\alpha(a) = \begin{bmatrix} a & 0 \\ 0 & a^{-1} \end{bmatrix} : a \in \mathbb{R}, a > 0\right\},$$
$$N = \left\{\beta(b) = \begin{bmatrix} 1 & b \\ 0 & 1 \end{bmatrix} : b \in \mathbb{R}\right\}.$$

Each $x \in G$ can be written uniquely as $x = \gamma(\theta)\alpha(a)\beta(b)$. Let $f \in UCS(G)$. For a fixed $\theta$, by replacing $f$ by $_{\gamma(-\theta)}f$, we may first consider $f \in UCS(AN)$. Note that $S = AN$ is a closed subgroup of $G$ and it is isomorphic to the group $G_3$ that we studied above. By Theorem 5.2, $f\big(\alpha(a)\beta(b)\big) = f(\alpha(a))$ for all $b \in \mathbb{R}$. By applying Lemma 4.3 and Lemma 4.4 of [5], one sees that $UCS(G)$ only contains constant functions. ∎

1616

Department of Mathematics, University at Buffalo, Buffalo, NY 14260-2900, USA

*Email address*: chouc@buffalo.edu